\documentclass{article}
\usepackage[utf8]{inputenc}
\usepackage{amsmath}\usepackage[hmargin=3cm,vmargin=3.5cm]{geometry}
\usepackage{graphicx}

\usepackage[dvipsnames]{xcolor}

\usepackage{tikz}

\usepackage{biblatex}
\usepackage{listings}

\definecolor{codegreen}{rgb}{0,0.6,0}
\definecolor{codegray}{rgb}{0.5,0.5,0.5}
\definecolor{codepurple}{rgb}{0.58,0,0.82}
\definecolor{backcolour}{rgb}{0.95,0.95,0.92}

\lstdefinestyle{mystyle}{
    backgroundcolor=\color{backcolour},   
    commentstyle=\color{codegreen},
    keywordstyle=\color{magenta},
    numberstyle=\tiny\color{codegray},
    stringstyle=\color{codepurple},
    basicstyle=\ttfamily\footnotesize,
    breakatwhitespace=false,         
    breaklines=true,                 
    captionpos=b,                    
    keepspaces=true,                 
    numbers=left,                    
    numbersep=5pt,                  
    showspaces=false,                
    showstringspaces=false,
    showtabs=false,                  
    tabsize=2
}

\title{Cantilevered, Rectangular Plate Dynamics \\ by Finite Difference Methods}
 \author{ 
            Benjamin Brown\\
 University of Maryland, Baltimore County \\
 \it Baltimore, MD \\
  \it bbrown4@umbc.edu}

\begin{document}

\maketitle

\abstract{In this technical note, we consider a dynamic linear, cantilevered rectangular plate. The evolutionary PDE model is given by the fourth order plate dynamics (via the spatial biharmonic operator) with clamped-free-free-free boundary conditions. We additionally consider damping/dissipation terms, as well as non-conservative lower order terms arising in various applications. Dynamical numerical simulations are achieved by way of a finite difference spatial approximation with a MATLAB time integrator. The rectangular geometry allows the use of standard 2D spatial finite differences, while the high spatial order of the problem and mixed clamped-free type boundary conditions present challenges. Dynamic energies are also computed. The relevant code is presented, with discussion of the model and context. }

\section{Introduction}

The purpose of this note is to describe a dynamic simulation method for a linear 2D cantilevered plate via finite differences. We specifically determine
the ghost points required for the enforcement of the cantilever boundary conditions (mixed clamped-free with corners). 

\subsection{Motivation}
The primary motivation for studying the particular PDE model below derives from the application of piezoelectric energy harvesting \cite{tang2019experimental, tang2018aeroelastic}. In this application, a cantilevered beam or rectangular plate is placed in a surrounding flow of air. In the {\em axial} configuration (whereby the flow runs from clamped end to free end of the structure), a cantilever is particularly prone to an aeroelastic instability termed {\em flutter}, even for low flow speeds. After the onset of instability, the beam settles into a limit cycle oscillation. If a piezo-electric device (for instance a thin layer, or patch) is affixed to the deflecting structure, power is generated and can be harvested. The energy is produced by the system's sustained
large, flapping motion which causes a current to be produced via the piezoelectric effect. Such energy-harvesting configurations have recently been shown to be viable \cite{tang2019experimental}. As an alternative energy source, scaled harvester systems could be effective in providing power to remote locations, or small surges of power to supplement a traditional grid. 

To effectively and efficiently harvest energy from the flow-induced oscillations of an elastic cantilever, one must have a viable distributed parameter system describing the large deflections of the 1D or 2D cantilevered structure. In particular, to produce a limit cycle oscillation in the post-flutter regime (i.e., after the onset of the flow-induced bifurcation), the structural model must permit large deflections and incorporate a nonlinear restoring force. Recent work has 1D {\em inextensible} cantilevers, including modeling \cite{dowell2016equations,deliyianni2020large}, and more recently, a well-posedness analysis of solutions \cite{deliyianni2021theory}. More recently, the modeling work was extended to (several) 2D inextensible cantilevered plates \cite{deliyianni2021dynamic}. These plate models are very involved, including nonlinear inertial and stiffness terms, as well as Lagrange multipliers for the enforcement of relevant constraints.  In this note, we focus on the linearization of the model in \cite{deliyianni2021dynamic}, and provide a fast numerical solver for a plate with given forcing (through the boundary and the right hand side). We accommodate damping effects, and crude approximations of the flow via a linear piston-theoretic approximation of potential flow (see, e.g., \cite{howell2019thorough}). This model is can be used viably in the pre-onset regime, as well as to compute a flow-induced bifurcation point via the non-conservative flow parameters. Additionally, the code is general enough to compute dynamic responses to arbitrary edge loading, as well as dynamic and spatially distributed right hand sides.

Much numerical work has been done on linear cantilevered models in 1D \cite{howell2018cantilevered,howell2019thorough,deliyianni2020large}, however much less
mathematical literature exists for 2D cantilevered plate models, likely owing to the challenges associated to the high order of the problem and requisitely mixed (with corners) boundary conditions\footnote{The mixed nature of the boundary conditions, as well as the corner angles, are known to affect the expected solution regularity---and also numerical convergence---near corners \cite{blum1980boundary}.}. We call attention to one older technical note that discusses finite difference stencils for rectangular plates with different boundary conditions \cite{barton1948finite}---the analysis there is entirely static. It is clear, in general, that  2D models (and associated simulations), are necessary for the engineering applications mentioned above. We resolve the PDE dynamics below using
 finite difference methods with careful calculations of ghost nodes that allow algebraic corner conditions to be properly solved. 

\subsection{Setup and Model}
An open rectangular domain, $\Omega$ with boundary $\Gamma=\overline{\Gamma_E} \cup \overline{\Gamma_N} \cup\overline{\Gamma_W} \cup\overline{\Gamma_S}$, will be used for this model, and will use cardinal coordinates for naming, as can be seen in the following diagram:
\begin{center}
\includegraphics[scale=.65]{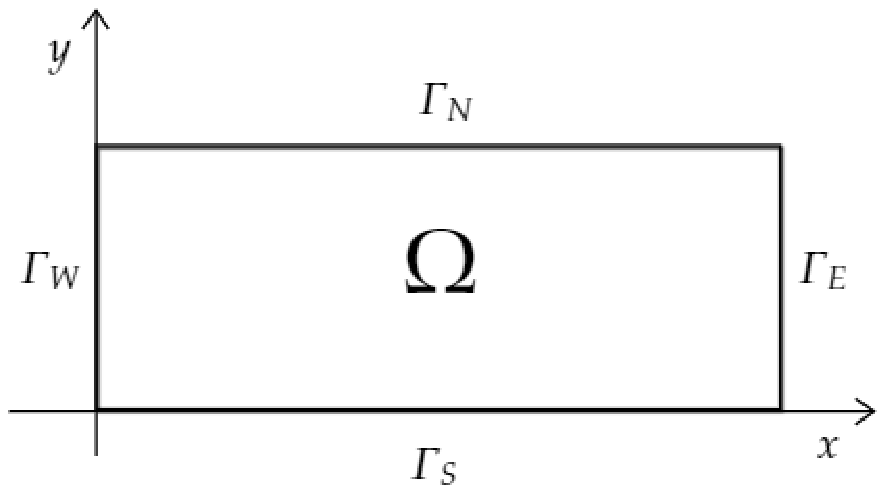}
\end{center}

The linear cantilevered plate model, with the rectangular domain $\Omega$, is described by the following equations:

\begin{equation}\label{PDE}
\begin{cases}
w_{tt} + D \Delta^2 w + k_0 w_t - k_1 \Delta w_t + a_1 w_x + a_2 w_y  = f(x,t)~~~&\text{in } \Omega  \\[0.5em]
w = \frac{\partial w}{\partial n} = 0 ~~~&\text{on }  \Gamma_W \\[0.5em]
\nu w_{xx}+w_{yy}=g_1(x,y),~~w_{yyy}+(2-\nu)w_{xxy}= h_1(x,y) ~~~&\text{on }  \Gamma_N   \\[0.5em]
\nu w_{xx}+w_{yy}=g_2(x,y),~~w_{yyy}+(2-\nu)w_{xxy}= h_2(x,y) ~~~&\text{on }  \Gamma_S  \\[0.5em]
w_{xx}+\nu w_{yy}=g_3(x,y),~~w_{xxx}+(2-\nu)w_{yyx}= h_3(x,y) ~~~&\text{on }  \Gamma_E \\[0.5em]
w(t=0)=w_0(x),~w_t(t=0)=w_1(x)
\end{cases}
\end{equation}
The quantity $\nu \in (0,1/2)$ represents the Poisson ratio \cite{deliyianni2021dynamic}.
Above, the function $w:\overline{\Omega} \times [0,T) \to \mathbf R$ represents the transverse displacement of an thin, isotropic, homogeneous plate. We do not explicitly incorporate piezoelectric effects here, and, assuming a piezoelectric material is attached to the surface of the plate, we incorporate the combined the combined stiffness effect through the parameter $D>0$.   The damping coefficients are $k_0,k_1\ge 0$, where the former measures weak or frictional damping, and the latter is so-called {\em square-root type damping} and is of a stronger nature. (See the discussion in \cite{howell2018cantilevered} for further discussion about the interpretation and physical meaning of this type of damping.) In reality, the term $-k_1\Delta w_t$ term may have questionable meaning for the cantilever, but is left in the implementation as it is unproblematic mathematically; we set $k_1=0$ for all  examples shown throughout this note.
The flow-velocity parameters are $a_1, a_2$ corresponding the respective $x$ and $y$ components of the surrounding airflow. A linear combination of these, for instance, represents
an off-axis flow; pure axial flow would take $a_1>0$ and $a_2=0$. The functions $g_i(x,y)$ correspond to edge-loading by way of boundary moments, while the boundary shears are given by $h_i(x,y)$ on each of the non-clamped boundaries in the model. The function $f(x,t)$ provides a distributed loading across the surface of $\Omega$.  In practice, for regular solutions, these boundary and interior data must satisfy natural compatibility conditions, in particular at/near corner points.

\begin{figure}
    \centering
    \begin{tikzpicture}[>=latex, scale=1.5]
        \draw[lightgray] (1,1) grid (5,5);
        
        \node[draw, minimum width=0.05cm, minimum height = 1cm,circle, fill=white, scale=0.75] (V1) at (3,3) {$O$};
        \node[draw, minimum width=0.05cm, minimum height = 1cm,circle, fill=white, scale=0.75] (V2) at (3,4) {$N$};
        \node[draw, minimum width=0.05cm, minimum height = 1cm,circle, fill=white, scale=0.75] (V3) at (3,5) {$NN$};
        \node[draw, minimum width=0.05cm, minimum height = 1cm,circle, fill=white, scale=0.75] (V4) at (3,2) {$S$};
        \node[draw, minimum width=0.05cm, minimum height = 1cm,circle, fill=white, scale=0.75] (V5) at (3,1) {$SS$};

        \node[draw, minimum width=0.05cm, minimum height = 1cm,circle, fill=white, scale=0.75] (V6) at (1,3) {$WW$};
        \node[draw, minimum width=0.05cm, minimum height = 1cm,circle, fill=white, scale=0.75] (V7) at (2,3) {$W$};
        \node[draw, minimum width=0.05cm, minimum height = 1cm,circle, fill=white, scale=0.75] (V8) at (4,3) {$E$};
        \node[draw, minimum width=0.05cm, minimum height = 1cm,circle, fill=white, scale=0.75] (V9) at (5,3) {$EE$};
        
        \node[draw, minimum width=0.05cm, minimum height = 1cm,circle, fill=white, scale=0.75] (V6) at (2,4) {$NW$};
        \node[draw, minimum width=0.05cm, minimum height = 1cm,circle, fill=white, scale=0.75] (V7) at (4,4) {$NE$};
        \node[draw, minimum width=0.05cm, minimum height = 1cm,circle, fill=white, scale=0.75] (V8) at (2,2) {$SW$};
        \node[draw, minimum width=0.05cm, minimum height = 1cm,circle, fill=white, scale=0.75] (V9) at (4,2) {$SE$};


    \end{tikzpicture}
    \caption{The 13-point stencil used in the finite difference method to approximate a fourth derivative}
    \label{fig:stencil}
\end{figure}
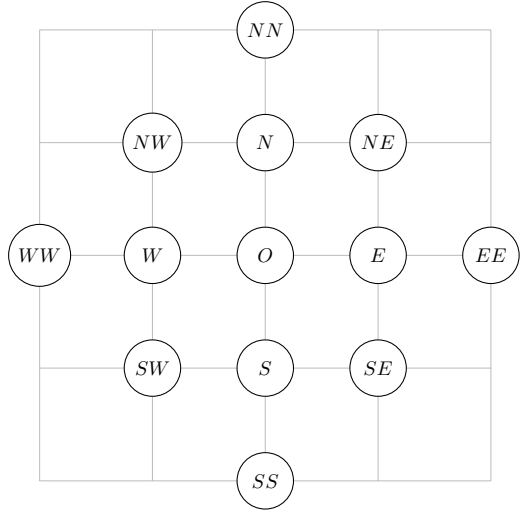

Due to the fourth order bi-harmonic operator, each point in the mesh created on the plate must be capable of having a 13-point stencil applied,
including the edges and corners, such as in Figure~\ref{fig:stencil}. Owing to this, ghost points must be added around the boundary of the plate.  To wit, a {\em ghost point} or {\em ghost node} is an artificial value
for a point that created outside of $\overline{\Omega}$ to enforce the boundary conditions, permitting a finite difference
stencil to be applied when the stencil overlaps boundary nodes. The particular challenges in the resolution and regularity of a given solution at corner points manifests themselves in algebraic issues for finite difference methods. Ghost points are calculated through using the boundary conditions, as well as known values from the interior, both of which are obtained through the data in the problem (e.g., $f$, $g_i$, and/or $h_i$). Indeed, the ghost node calculations are a central contribution of this note. 

\section{Ghost Points}

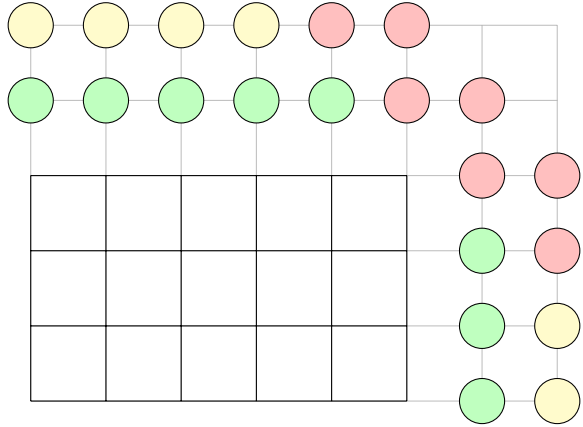
\begin{figure}
    \centering
    \begin{tikzpicture}[>=latex, scale=1]
        \draw[lightgray] (1,1) grid (8,6);
        \draw[black] (1,1) grid (6,4);
        
        \draw[fill = green!25] (1,5) circle (0.3);
        \draw[fill = green!25] (2,5) circle (0.3);
        \draw[fill = green!25] (3,5) circle (0.3);
        \draw[fill = green!25] (4,5) circle (0.3);
        \draw[fill = green!25] (5,5) circle (0.3);
        
        \draw[fill = green!25] (7,3) circle (0.3);
        \draw[fill = green!25] (7,2) circle (0.3);
        \draw[fill = green!25] (7,1) circle (0.3);
        
        \draw[fill = yellow!25] (1,6) circle (0.3);
        \draw[fill = yellow!25] (2,6) circle (0.3);
        \draw[fill = yellow!25] (3,6) circle (0.3);
        \draw[fill = yellow!25] (4,6) circle (0.3);

        \draw[fill = yellow!25] (8,2) circle (0.3);
        \draw[fill = yellow!25] (8,1) circle (0.3);
        
        \draw[fill = red!25] (5,6) circle (0.3);
        \draw[fill = red!25] (6,5) circle (0.3);
        \draw[fill = red!25] (6,6) circle (0.3);
        \draw[fill = red!25] (7,5) circle (0.3);
        \draw[fill = red!25] (7,4) circle (0.3);
        \draw[fill = red!25] (8,4) circle (0.3);
        \draw[fill = red!25] (8,3) circle (0.3);

    \end{tikzpicture}
    \caption{Diagram showing the order of ghost points being calculated around a free-free corner. Begin by calculating the first row, which are the green
    nodes.  Then the second row will be calculated, represented by yellow nodes. Finally the system of seven equations can be solved
    for the corner nodes, shown in red.}
    \label{fig:nodeorder}
\end{figure}

As alluded to above, the inherent mixed boundary conditions in the cantilever problem introduce complications into the calculation of the ghost points. To determine ghost node values, we will begin with the clamped edge, because it is simple and only has one unknown. This will provide information to supply to more stencils in latter calculations of nodal values.

The nomenclature used in this section makes use of {\em cardinal values} that correspond to the current displacement of the point in the stated direction,
with $O$ representing the `origin'. For example, $E$ refers to the mesh point to the right of $O$.

\subsection{Clamped Edge}

On this edge, the conditions require that $w = 0$ and $w_x = 0$.  Using the finite difference approximation on this point:
\[
E = W,
\]
for all $O$ along the clamped boundary.  

\subsection{Free Edge}
All of the following calculations will be made for the $\Gamma_N$ boundary, but they apply to the other edges
with a simple rotation of coordinates as needed.

The conditions on the three free edges require specified moments and shears:
\[
\begin{cases}
\nu w_{xx}+w_{yy}=g_1(x,y),~~w_{yyy}+(2-\nu)w_{xxy}= h_1(x,y) ~~~&\text{on }  \Gamma_N   \\[0.5em]
\nu w_{xx}+w_{yy}=g_2(x,y),~~w_{yyy}+(2-\nu)w_{xxy}= h_2(x,y) ~~~&\text{on }  \Gamma_S  \\[0.5em]
w_{xx}+\nu w_{yy}=g_3(x,y),~~w_{xxx}+(2-\nu)w_{yyx}= h_3(x,y) ~~~&\text{on }  \Gamma_E \\[0.5em]
\end{cases}
\]
Because a 13-point stencil will be needed for the 
final calculation, two rows of ghost points must be determined.

\subsubsection{First Row on Edge: Avoiding Corners}
The second order boundary condition will be used for the first row.  The
free-free corner will be avoided for this calculation, as it has has special calculations, shown later.  Using
the finite difference method to approximate  this condition yields the following:
expression:
\[
\label{eq:bc1}
{\frac {\nu\, \left( W-2\,O+E \right) }{{{dx}}^{2}}}+{\frac {N-2\,
O+S}{{{dy}}^{2}}}= g_1(x,y).  
\]

Notice that along all the free edges, except when adjacent to a free-free corner, there is only one 
unknown value, the ghost value.  This can be explicitly solved for, for example on the southern edge, to
give the requisite value of the ghost point:
\[
S=-{\frac {\nu\, \left( W-2\,O+E \right) {{dy}}^{2}}{{{dx}}^{2
}}}-N+2\,O + g_1(x,y).
\]

\subsubsection{Second Row on Edge: Avoiding Corners}
The second boundary condition will be used for the second row of ghost points.  When this stencil is utilized, and the
first edge has already been calculated, the second row ghost values can then be solved directly.  This calculations uses
ghost nodes found in the previous step in the first row.  Because of this, there must again be caution about the corner.
The nodes to be calculated are shown in Figure~\ref{fig:nodeorder} as the yellow
nodes. Using the finite difference approximation on the second boundary condition yields the following expression:
\[
{\frac {{SS}-2\,S+2\,N-{NN}}{2\,{{dy}}^{3}}} + 
    {\frac {\left( 2-\nu \right)  \left( {SE}-{NE}-2\,S+2\,N+{SW}-{NW} \right) }{2\,{{dx}}^{2}{dy}}}= h_1(x,y).
\]
The ghost point $SS$ can be determined  from the above, as it is the only unknown when applied to the southern free edge.  The expression
for this value is:
\[
{SS}=-{\frac { \left( 2-\nu \right) {{dy}}^{2} \left( {SE}
-{NE}-2\,S+2\,N+{SW}-{NW} \right) }{{{dx}}^{2}}}+2\,S-
2\,N+{NN} + h_1(x,y).
\]

\subsubsection{Corners}

For each of the two free-free corners, there are seven unknown ghost points.  We can make a linear system of the seven boundary conditions that
can be applied to the surrounding points in order to explicitly find each of them.  In this section, {\bf for simplicity, the 
$g_i(x,y)$ and $h_i(x,y)$ terms will be set to zero}.  In a latter section we will elaborate upon the edge loading conditions.
 
The relevant finite difference boundary conditions are:
\[\label{corner}
\begin{cases}
{\frac {\nu\, \left( W-2\,O+E \right) }{{{dx}}^{2}}}+{\frac {N-2\,O+S}{{{dy}}^{2}}}=0  \\[0.5em]
{\frac {W-2\,O+E}{{{dx}}^{2}}}+{\frac {\nu\, \left( N-2\,O+S\right) }{{{dy}}^{2}}}=0 \\[0.5em]
{\frac {{SS}-2\,S+2\,N-{NN}}{2\,{{dy}}^{3}}}+{\frac {\left( 2-\nu \right)  \left( {SE}-
{NE}-2\,S+2\,N+{SW}-{NW} \right) }{2\,{{dx}}^{2}{dy}}}=0 \\[0.5em]
{\frac {{SSW}-2\,{SW}+2\,{NW}-{NNW}}{2\,{{dy}}^{3}}}+{\frac { \left( 2-\nu \right)  
\left( S-N-2\,{SW}+2\,{NW}+{SWW}-{NWW} \right) }{2\,{{dx}}^{2}{dy}}}=0 \\[0.5em]
{\frac {{EE}-2\,E+2\,W-{WW}}{2\,{{dx}}^{3}}}+{\frac {\left( 2-\nu \right)  \left( {SE}-
{SW}-2\,E+2\,W+{NE}-{NW} \right) }{2\,{{dy}}^{2}{dx}}}=0 \\[0.5em]
{\frac {{SEE}-2\,{SE}+2\,{SW}-{SWW}}{2\,{{dx}}^{3}}}+{\frac { \left( 2-\nu \right)  
\left( {SSE}-{SSW}-2\,{SE}+2\,{SW}+E-W \right) }{2\,{{dy}}^{2}{dx}}}=0 \\[0.5em]
{SE}-{NE}-{SW}+{NW}=0 \\[0.5em]
\end{cases}
\]
Solving this system of boundary conditions yields the expressions for each ghost point about corner $O$,
when $O$ is the corner between $\Gamma_N$ and $\Gamma_E$:
\[
\begin{cases}
E=2\,O-W  \\[0.5em]
N=2\,O-S\\[0.5em]
{NE}={SE}-{SW}+{NW}\\[0.5em]
EE= \frac{(-4O+4W+2SE-2SW )(\nu-2){dx}^2}{{dy}^2} + ( 4O-4W+WW ) \\[0.5em]
NN= \frac{(-4O+4S+2NW-2SW )(\nu-2){dy}^2}{{dx}^2} + ( 4O-4S+SS ) \\[0.5em]
\it{NNW}=\frac {( \nu-2 )  ( 2O-SWW-2S-2NW+2SW+NWW ) {dy}^2}{{dx}^2} + ( 2NW-2SW+SSW) \\[0.5em]
\it{SEE}=\frac {( \nu-2 )  ( 2O-W-2SE+2SW+2SSE-SSW ) {dy}^2}{{dx}^2} + ( 2SWW+2SE-2SW) \\[0.5em]
\end{cases}
\]
The same approach can be applied to the $\Gamma_S$ and $\Gamma_E$ corner, and the results are simply
a reflection of these expressions.

\section{Non-zero Edge Conditions}

In the previous section the $g_i(x,y)$ and $h_i(x,y)$ values were set to be zero.  This was done for simplicity, and
clarity of exposition; in this section we provide examples of ghost node calculations with constant edge loading. As in the theory of boundary lifts, we can compute the solution from a given edge loading, while setting all other boundary loads to zero. The response to several edge loads can thus be reconstructed for the linear problem through the principle of superposition.

  In this section, a derivation of nonzero, constant values will be shown, and the result of these loads will be determined.  Again we shall only  solve
for the $\Gamma_N$ boundary, but other boundary values will be determined via rotation.

We first handle the case where $g_i(x,y) = G$ and $h_i(x,y) = 0$.  The conditions for this case become:
\[\label{cornerG}
\begin{cases}
{\frac {\nu\, \left( W-2\,O+E \right) }{{{dx}}^{2}}}+{\frac {N-2\,O+S}{{{dy}}^{2}}}=G  \\[0.5em]
{\frac {W-2\,O+E}{{{dx}}^{2}}}+{\frac {\nu\, \left( N-2\,O+S\right) }{{{dy}}^{2}}}=G \\[0.5em]
{\frac {{SS}-2\,S+2\,N-{NN}}{2\,{{dy}}^{3}}}+{\frac {\left( 2-\nu \right)  \left( {SE}-
{NE}-2\,S+2\,N+{SW}-{NW} \right) }{2\,{{dx}}^{2}{dy}}}=0 \\[0.5em]
{\frac {{SSW}-2\,{SW}+2\,{NW}-{NNW}}{2\,{{dy}}^{3}}}+{\frac { \left( 2-\nu \right)  
\left( S-N-2\,{SW}+2\,{NW}+{SWW}-{NWW} \right) }{2\,{{dx}}^{2}{dy}}}=0 \\[0.5em]
{\frac {{EE}-2\,E+2\,W-{WW}}{2\,{{dx}}^{3}}}+{\frac {\left( 2-\nu \right)  \left( {SE}-
{SW}-2\,E+2\,W+{NE}-{NW} \right) }{2\,{{dy}}^{2}{dx}}}=0 \\[0.5em]
{\frac {{SEE}-2\,{SE}+2\,{SW}-{SWW}}{2\,{{dx}}^{3}}}+{\frac { \left( 2-\nu \right)  
\left( {SSE}-{SSW}-2\,{SE}+2\,{SW}+E-W \right) }{2\,{{dy}}^{2}{dx}}}=0 \\[0.5em]
{SE}-{NE}-{SW}+{NW}=0 \\[0.5em]
\end{cases}
\]
\noindent Solving this system for the conditions, as done previously, provides the solution:
\[
\begin{cases}
E=2\,O-W + \frac{G dx^2}{\nu+1} \\[0.5em]
N=2\,O-S+\frac{G dy^2}{\nu+1}\\[0.5em]
NE=SE-SW+NW\\[0.5em]
EE= \frac{(-4O+4W+2SE-2SW )(\nu-2){dx}^2}{{dy}^2} + ( 4O-4W+WW ) 
- {\frac {2G{{dx}}^{4}\nu - 4G{{dx}}^{4}}{{{dy}}^{2} \left( \nu+1 \right) }
}+{\frac {2G{{dx}}^{2}}{\nu+1}}\\[0.5em]
NN= \frac{(-4O+4S+2NW-2SW )(\nu-2){dy}^2}{{dx}^2} + ( 4O-4S+SS )
- {\frac {2G{{dy}}^{4}\nu - 4G{{dy}}^{4}}{{{dx}}^{2} \left( \nu+1 \right) }
}+{\frac {2G{{dy}}^{2}}{\nu+1}}\\[0.5em]
NNW=\frac {( \nu-2 )  ( 2O-SWW-2S-2NW+2SW+NWW ) {dy}^2}{{dx}^2} + ( 2NW-2SW+SSW)
+{\frac {G \left( -2+\nu \right) {{dy}}^{4}}{{{dx}}^{2}
 \left( \nu+1 \right) }}\\[0.5em]
SEE=\frac {( \nu-2 )  ( 2O-W-2SE+2SW+2SSE-SSW ) {dy}^2}{{dx}^2} + ( 2SWW+2SE-2SW)
+{\frac {G \left( -2+\nu \right) {{dx}}^{4}}{{{dy}}^{2}
 \left( \nu+1 \right) }}
\\[0.5em]
\end{cases}
\]
Similarly we can consider the case where $g_i(x,y) = 0$ and $h_i(x,y) = H$.  The conditions
for this system become:
\[\label{cornerH}
\begin{cases}
{\frac {\nu\, \left( W-2\,O+E \right) }{{{dx}}^{2}}}+{\frac {N-2\,O+S}{{{dy}}^{2}}}=0  \\[0.5em]
{\frac {W-2\,O+E}{{{dx}}^{2}}}+{\frac {\nu\, \left( N-2\,O+S\right) }{{{dy}}^{2}}}=0 \\[0.5em]
{\frac {{SS}-2\,S+2\,N-{NN}}{2\,{{dy}}^{3}}}+{\frac {\left( 2-\nu \right)  \left( {SE}-
{NE}-2\,S+2\,N+{SW}-{NW} \right) }{2\,{{dx}}^{2}{dy}}}=H \\[0.5em]
{\frac {{SSW}-2\,{SW}+2\,{NW}-{NNW}}{2\,{{dy}}^{3}}}+{\frac { \left( 2-\nu \right)  
\left( S-N-2\,{SW}+2\,{NW}+{SWW}-{NWW} \right) }{2\,{{dx}}^{2}{dy}}}=H \\[0.5em]
{\frac {{EE}-2\,E+2\,W-{WW}}{2\,{{dx}}^{3}}}+{\frac {\left( 2-\nu \right)  \left( {SE}-
{SW}-2\,E+2\,W+{NE}-{NW} \right) }{2\,{{dy}}^{2}{dx}}}=H \\[0.5em]
{\frac {{SEE}-2\,{SE}+2\,{SW}-{SWW}}{2\,{{dx}}^{3}}}+{\frac { \left( 2-\nu \right)  
\left( {SSE}-{SSW}-2\,{SE}+2\,{SW}+E-W \right) }{2\,{{dy}}^{2}{dx}}}=H \\[0.5em]
{SE}-{NE}-{SW}+{NW}=0 \\[0.5em]
\end{cases}
\]
\noindent The solutions for this system are:
\[
\begin{cases}
E=2\,O-W  \\[0.5em]
N=2\,O-S\\[0.5em]
NE={SE}-{SW}+{NW}\\[0.5em]
EE= \frac{(-4O+4W+2SE-2SW )(\nu-2){dx}^2}{{dy}^2} + ( 4O-4W+WW ) + 2Hdx^3 \\[0.5em]
NN= \frac{(-4O+4S+2NW-2SW )(\nu-2){dy}^2}{{dx}^2} + ( 4O-4S+SS )  - 2Hdx^3 \\[0.5em]
NNW=\frac {( \nu-2 )  ( 2O-SWW-2S-2NW+2SW+NWW ) {dy}^2}{{dx}^2} + ( 2NW-2SW+SSW) - 2Hdx^3 \\[0.5em]
SEE=\frac {( \nu-2 )  ( 2O-W-2SE+2SW+2SSE-SSW ) {dy}^2}{{dx}^2} + ( 2SWW+2SE-2SW)  + 2Hdx^3 \\[0.5em]
\end{cases}
\]

\section{Reduction of Order and ODE Solver}
To integrate in time, and produce distributed values at each time step, we invoke
 MATLAB's {\tt ode15s} solver.  This solver requires an abstract differential equation of the form:
\begin{equation}
    \dot y= F(t,y).
\end{equation}

To reduce our evolutionary PDE, we introduce states to reduce the order and create a $2\times 2$ system in the state
variables to obtain a first order formulation.  The problem also has many spatial nodes, which each must be resolved. 
A matrix expression will be used to obtain one principal variable $y$ to feed to the solver.  We begin with the interior PDE:
\begin{equation}
w_{tt} + D \Delta^2 w + a_1 w_x + a_2 w_y + k_0 w_t - k_1 (w_{txx} + w_{tyy}) = f(x,t),
\end{equation}
and introduce the state variable $v$, capturing $w_t$, or the velocity.  We can also
take the derivative of $v$ and using the PDE expression for $w_{tt}$, we obtain the value of $v_t$:
\begin{equation}\begin{cases}
    v =&w_t\\
    v_t =&  f(x,t) - D \Delta^2 w - a_1 w_x - a_2 w_y - k_0 w_t + k_1 (w_{txx} + w_{tyy}).
    \end{cases}
\end{equation}
We can write each of these as a matrix:
\[
y = \begin{bmatrix} w \\ v \end{bmatrix}, y_t = \begin{bmatrix} w_t \\ v_t \end{bmatrix}
\]
For each node in the model we must solve for $w$ and $v$, and thus we will extend each of these in the 
array with indices corresponding to each point in the mesh. There are a total of $N$ nodes, determined 
from user set parameters as laid out in Section~\ref{sec:variables}.
\[
y = \begin{bmatrix} w_1 , w_2 , \hdots , w_N , v_1 , v_2 , \hdots, v_N \end{bmatrix}^T
\]
Using these vectors, this problem is now written as a first order ODE, and can be solved directly with the {\tt ode15s} solver,
where each iteration will be stored for future animation, post-processing calculations, and other data analyses.  Appendix \ref{sec:appendix} contains MATLAB code used to simulate the plate dynamics. 

\section{Energy Calculation}
The natural potential energy associated to a free plate is expressed in terms of the bilinear form
\begin{equation}\label{weak}
a(u,v) =  \nu ( \Delta u, \Delta v) + (1-\nu) \big[( u_{xx}, v_{xx}) + 2(u_{xy},v_{xy}) + (u_{yy}, v_{yy})  \big],
\end{equation}
where the expression $(f,g)$ represents the $L^2(\Omega)$ inner-product:
$$(f,g)=(f,g)_{L^2(\Omega)}=\int_{\Omega} f(x,y)g(x,y)d\mathbf x.$$

We consider the standard kinetic and potential contributions to the total energy. This is to say, for a state $[w,v]$ we can compute the total energy as 
$$E(t)=E(w(t),v(t))=U(w(t))+K(v(t)).$$

The potential energy is given through the bilinear form as $U(w)=a(w,w)$:
\begin{equation*}
U(w(t)) = \frac{1}{2} D  \int_{\Omega} \left[ \nu(\Delta w(t))^2 +(1-\nu) \left( w^2_{xx}(t) + 2w(t)^2_{xy} + w(t)^2_{yy}\right) \right] d\mathbf x 
\end{equation*}
and the kinetic energy 
\begin{equation*}
K(v) = \frac{1}{2}   \int_{\Omega} v(t)^2 d\mathbf x,
\end{equation*}
where we have suppressed the dependence of the states on $\mathbf x$.

\begin{figure}
    \centering
    \includegraphics[width=\linewidth]{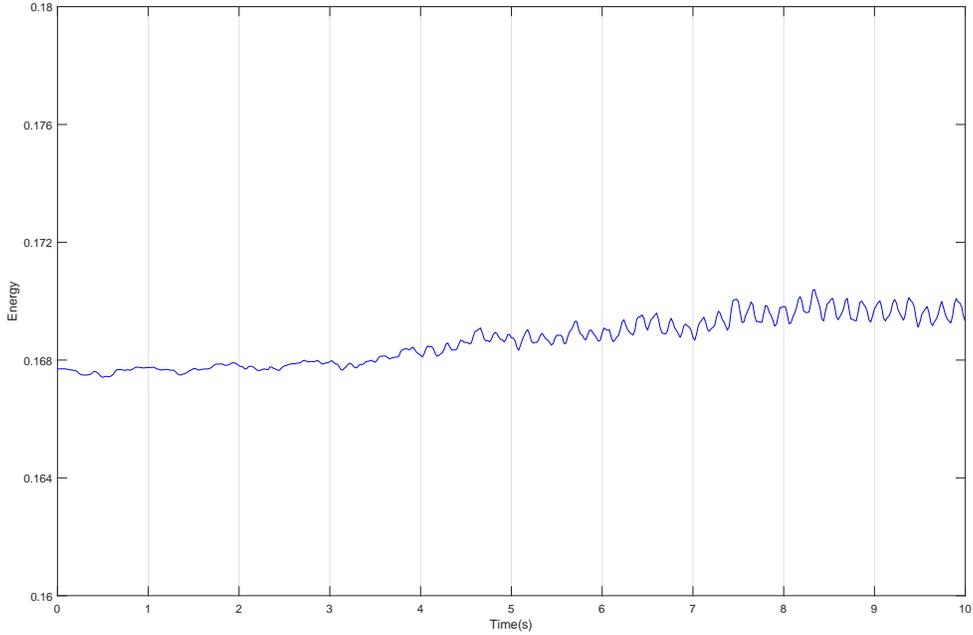}
    \caption{Conservation of energy demonstrated with a 10x10 mesh grid and $k_0=k_1=a_1=a_2=0$, $v(x,y;0)=x$}
    \label{fig:consEnergy}
\end{figure}

If there is no damping or flow considered in the model (i.e., $k_0=k_1=a_1=a_2=0$), then conservation of energy is expected for solutions---see, for instance, \cite{deliyianni2021theory,howell2019thorough}. We demonstrate this here with an initial state consisting of an initial displacement of $w(x,y;0)=0$ and an initial velocity of $v(x,y;0)=x$.
Figure~\ref{fig:consEnergy} is produced by the script and demonstrates the conservation of energy for the numerically computed solution in this situation.

Of course, slight variations (errors) can be observed on this small scale due to numerical error in the finite difference method and time integrator used.

\begin{figure}
    \centering
    \includegraphics[width=\linewidth]{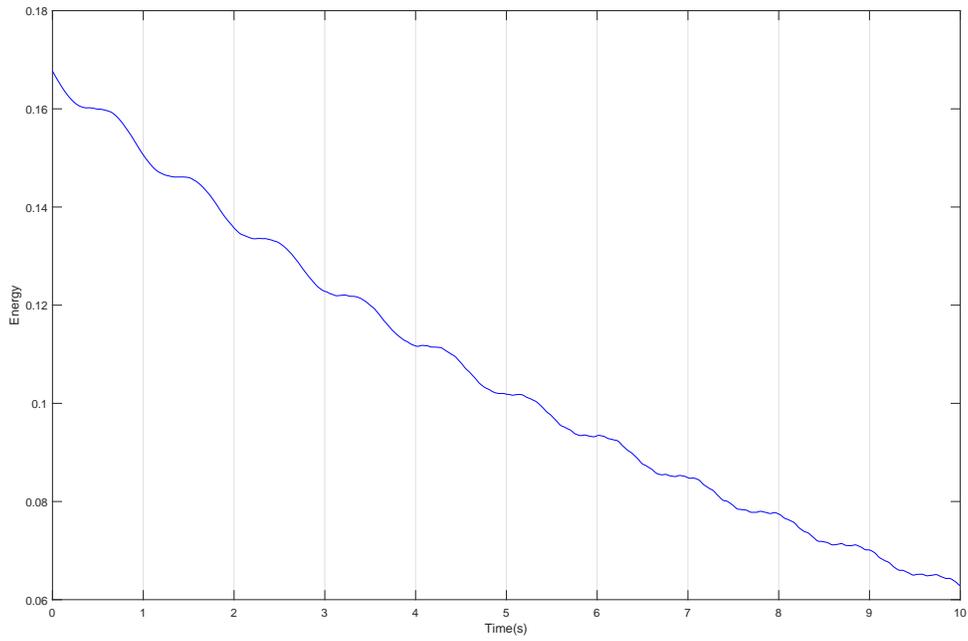}
    \caption{Decay of energy demonstrated with a 10x10 mesh grid and $k_1=a_1=a_2=0$, $k_0=0.1$ and $v(x,y;0)=x$}
    \label{fig:dampEnergy}
\end{figure}

We also demonstrate decay of the energy in the presence of  some damping.  In a similar case to the above, but taking a small amount of frictional damping $k_0=0.1$, we observe the energy decreasing.  Figure~\ref{fig:dampEnergy} is produced by the script and demonstrates the dissipation of energy in the system.

\newpage

\section{Table of Variables}
\label{sec:variables}
In the script that was generated in correspondence to this document, there are many variables that can 
be modified in order to fit the desired situation.  Table \ref{tab:Variables} contains a list of the names 
and descriptions of each of these variables. 

\begin{table}[ht]
    \centering
    \begin{tabular}{c|c}
         Variable & Description\\
         \hline
         D & Stiffness parameter for the cantilever material\\
         Lx & Length of the cantilever in the x direction\\
         Ly & Length of the cantilever in the y direction\\
         Nx & Number of points in the mesh in the x direction\\
         Ny & Number of points in the mesh in the y direction\\
         N & Total number of node points where w is unknown\\
         $\nu$ & Poisson Ratio\\
         $k_0$ & Linear damping parameter\\
         $k_1$ & Nonlinear damping parameter\\
         $a_1$ & Flow parameter in the x direction\\
         $a_2$ & flow parameter in the y direction\\
         anim & Turns the animation on/off\\
         energies & Turns the energy calculation on/off\\
         $t_0$ & Initial time of the calculation\\
         $t_f$ & Final time of the calculation\\
         ns & Number of time steps taken\\
         
    \end{tabular}
    \caption{Descriptions of each of the variables that can be chosen in the code}
    \label{tab:Variables}
\end{table}

In addition to these variables, the user can also set the desired initial conditions using the functions described in Table \ref{tab:Functions}.

\begin{table}[ht]
    \centering
    \begin{tabular}{c|p{5cm}}
         Function & Description\\
         \hline
         winit(X, Y) & Sets the initial position of the cantilever.  This is done by defining the internal variable
                        wmat as a function of the X and Y positional arrays.\\
         vinit(X, Y) & Sets the initial velocity of the cantilever.  This is done by defining the internal variable
                        vmat as a function of the X and Y positional arrays.\\

    \end{tabular}
    \caption{Descriptions of each of the initial condition functions that can be chosen in the code}
    \label{tab:Functions}
\end{table}

\section{Acknowledgements}
The author acknowledges the generous support of the National Science Foundation in this work through the grant entitled: ``Collaborative Research: Experiment, Theory, and Simulation of Aeroelastic Limit Cycle Oscillations for Energy Harvesting Applications" (NSF-DMS 1908033).
\pagebreak

\appendix

\section{Matlab Code}
\label{sec:appendix}
\begin{lstlisting}[language=Matlab]
close all;

%problem parameters
D = 1.0e0; %stiffness parameter
Lx = 1; %length in x direction
Ly = 1; %length in y direction
Nx = 10; %number of nodes in x direction
Ny = 10; %number of nodes in y direction
N = (Nx-1)*(Ny); %total number of nodes where w is unknown
nu=0.3;
k0 = 0.1; %damping parameter
k1 = 0; %damping parameter
a1 = 0;
a2 = 0;
anim = 1; %see animation
energies = 1; %compute and plot energies

%timestepping parameters
t0 = 0; %initial time
tf = 10;
0e0; %final time
ns = ceil(50*tf)+1; %number of time steps
t = linspace(t0, tf, ns); %time vector


%initialize spatial domain
x = linspace(0,Lx,Nx);
y = linspace(0,Ly,Ny);
[X,Y] = meshgrid(x,y);
dx = x(2)-x(1);
dy = y(2)-y(1);

%initialize solution array W = [wn; vn];
%get initial diplacement and velocity - form W0
wn = winit(X,Y);
vn = vinit(X,Y);
W0 = [wn; vn]; %initial data vector is of length 2*N


% dydt = RHS(t,W0,D,Nx,Ny,N,x,y,dx,dy,nu);
% dydt = RHS(t,dydt,D,Nx,Ny,N,x,y,dx,dy,nu);
% dydt = RHS(t,dydt,D,Nx,Ny,N,x,y,dx,dy,nu);
% dydt = RHS(t,dydt,D,Nx,Ny,N,x,y,dx,dy,nu);

% call the timestepping code
tic;
[T,w] = ode15s(@(t,W) RHS(t,W,D,Nx,Ny,N,x,y,dx,dy,nu,k0,k1,a1,a2),t,W0);
odetime = toc;
fprintf('Total run time for ODE integrator: %g\n',odetime);

%postprocess the results
%plot the surface
wmat = [zeros(Ny,1) reshape(w(end,1:N),[Nx-1 Ny])'];

% figure(3); surf(X,Y,wmat); %plot the final surface
% zlabel({'w(x,y)'});
% ylabel({'y'});
% xlabel({'x'});
% title('Solution at final time');

%animate the solution
wvals=w(:,1:N);
wmax=max(max(abs(wvals)));
if anim>0
    fig = figure(11);%hold on
    F(ns) = struct('cdata',[],'colormap',[]);
    grid on;
    plt=surf(X,Y,[zeros(Ny,1) reshape(w(1,1:N),[Nx-1 Ny])' ]);
    colormap(jet(256));
    colorbar;
    caxis([-1.1*wmax,1.1*wmax]);
    zlim([-1.1*wmax,1.1*wmax])
    view(-40,30);
    str=sprintf('t = %3.3f',T(1));
    h=text(0.95,0.95,0.9*wmax,str,'FontSize',12);
    F(1) = getframe(fig);
    for j = 2:ns
        plt.ZData =[zeros(Ny,1) reshape(w(j,1:N),[Nx-1 Ny])'];
        str=sprintf('t = %3.3f',T(j));
        set(h,'String',str);
        drawnow % display updates
        F(j) = getframe(fig);
    end
    movie(fig,F,2,anim*ns/tf);
    
    vid = VideoWriter('animation.mp4','MPEG-4');
    open(vid);
    writeVideo(vid,F);
    close(vid);
end

% compute energies if desired
if energies==1
    tic;
    %compute and plot energies
    Le = zeros(size(ns));
    Nle = zeros(size(ns));
    for j=1:ns
        Le(j)= computeEnergies(w(j,:),Nx,Ny,x,y,dx,dy,D,nu);
    end
    energytime = toc;
    fprintf('Total run time for energy calculation: %g\n',energytime);
    figure;hold
    plot(T,Le,'-b', 'DisplayName','E(t)');
    legend('linear E(t)')
end

%end



%function for initial displacement
function w = winit(X,Y)

Lx = X(1,end); %get Lx
Ly = Y(end,1); %get Ly
[Ny,Nx] = size(X); %get Nx and Ny
%wmat = (-4*X.^5+15*X.^4-20*X.^3+10*X.^2); %define the initial value of w
wmat = 0.0*X; %define the initial value of w
% figure; surf(X,Y,wmat); %plot the initial surface
% zlabel({'w(x,y)'});
% ylabel({'y'});
% xlabel({'x'});
% title('Initial Displacement');
wmat = wmat(:,2:Nx); %trim off left values
w = reshape(wmat',[],1); %reshape matrix into the w vector
end

%function for initial velocity
function v = vinit(X,Y)
Lx = X(1,end); %get Lx
Ly = Y(end,1); %get Ly
[Ny,Nx] = size(X); %get Nx and Ny
vmat = X; %define the initial value of w
% figure; surf(X,Y,vmat); %plot the initial surface
% % Create zlabel
% zlabel({'v(x,y)'});
% % Create ylabel
% ylabel({'y'});
% % Create xlabel
% xlabel({'x'});
% title('Initial Velocity');
vmat = vmat(:,2:Nx); %trim off left
v = reshape(vmat',[],1); %reshape matrix into the w vector
end

%function for computing energies
function [le,nle] = computeEnergies(w,Nx,Ny,x,y,dx,dy,D,nu)

%separate the displacements from velocities
w = reshape(w,[Nx-1 2*(Ny)])';
wmat = [zeros(Ny,1),w(1:Ny,:)]; %displacements at interior nodes
vmat = [zeros(Ny,1),w(Ny+1:2*Ny,:)]; %velocities at interior nodes

%compute the appropriate derivatives at the interior points
[lx,ly,lxy] = laplace(wmat,Nx,Ny,dx,dy,nu);
%compute the integral of grad w squared
integrand = nu * (lx + ly).^2 + (1-nu)*(lx.^2 + 2.*lxy.^2 + ly.^2);
normpotential2 = trapz(y,trapz(x,integrand,2));

%compute the integral of the velocities squared
normvel2 = trapz(y,trapz(x,vmat.*vmat,2));

%compute linear energy
le = 0.5*D*normpotential2 + 0.5*normvel2;

end
function [lx,ly,lxy] = laplace(w,Nx,Ny,dx,dy,nu)

%initialize the output arrays
lx = zeros(size(w));
ly = zeros(size(w));
lxy = zeros(size(w));
%apply BCs to wmat so we can compute w_xx and w_yy
%build the array with ghost nodes
w = [zeros(1,size(w,2)+2); [zeros(size(w,1),1), w, zeros(size(w,1),1)]; zeros(1,size(w,2)+2)];

%Clamped edge ghost points
w(:,1) = w(:,3);

%free edge (avoiding free-free corner area)
for j=2:Nx
    i = Ny+1;
    w(i+1,j) = -nu * (w(i,j-1) - 2 * w(i,j) + w(i,j+1)) / dx ^ 2 * dy ^ 2 - w(i-1,j) + 2 * w(i,j);
end

for j=2:Nx
    i = 2;
    w(i-1,j) = -nu * (w(i,j-1) - 2 * w(i,j) + w(i,j+1)) / dx ^ 2 * dy ^ 2 - w(i+1,j) + 2 * w(i,j);
end

for i=3:Ny+1
    j = Nx+1;
    w(i,j+1) = -nu * (w(i-1,j) - 2 * w(i,j) + w(i+1,j)) / dy ^ 2 * dx ^ 2 - w(i,j-1) + 2 * w(i,j);
end


%Now we will apply the conditions for the corners.
%Labelling convention in this section is relative to the corner

%Top right corner
i = 2;
j = Nx + 1;
%E
w(i,j+1) = 2 * w(i,j) - w(i,j-1);

%N
w(i-1,j) = 2 * w(i,j) - w(i+1,j);

%NE
w(i-1,j+1) = w(i+1,j+1) - w(i+1,j-1) + w(i-1,j-1);


%Bottom right corner
i = Ny + 1;
j = Nx + 1;
%E
w(i,j+1) = 2 * w(i,j) - w(i,j-1);


%S
w(i+1,j) = 2 * w(i,j) - w(i-1,j);

%SE
w(i+1,j+1) = w(i-1,j+1) - w(i-1,j-1) + w(i+1,j-1);

%loop over interior nodes and compute w_x and w_y using 2nd order centered
%differences
for j = 2:Nx+1 %which row of nodes
    for i=2:Ny+1 %which column of nodes
        lx(i-1,j-1) = (1/(dx^2))*(w(i,j+1)-2*w(i,j)+w(i,j-1)); %w_xx calculation
        ly(i-1,j-1) = (1/(dy^2))*(w(i+1,j)-2*w(i,j)+w(i-1,j)); %w_yy calculation
        lxy(i-1,j-1) = (1/(4*dx*dy)) * ( w(i+1,j+1) - w(i-1,j+1) - w(i+1,j-1) + w(i-1,j-1));
    end
end
end

function [gx,gy] = grads(w,Nx,Ny,dx,dy)

gx = zeros(size(w));
gy = zeros(size(w));

%apply BCs to wmat so we can compute w_xx and w_yy
%build the array with ghost nodes
w = [zeros(1,size(w,2)+2); [zeros(size(w,1),1), w, zeros(size(w,1),1)]; zeros(1,size(w,2)+2)];
%Applying hinged-hinged conditions
w(1,:) = -1.0 * w(3,:);
w(end,:) = -1.0 * w(end-2,:);

w(:,1) = -1.0 * w(:,3);
w(:,end) = -1.0 * w(:,end-2);
%loop over interior nodes and compute w_x and w_y using 2nd order centered
%differences
for j = 2:Ny+1 %which row of nodes
    for i=2:Nx-1 %which column of nodes
        gx(j-1,i) = (1/(2*dx))*(w(j,i+1)-w(j,i-1)); %w_x calculation
        gy(j-1,i) = (1/(2*dy))*(w(j+1,i)-w(j-1,i)); %w_y calculation
    end
end
end

%define the forcing function f(x,y,t)
function z = f(~,~,~)
z=0;
end

%define the ODE RHS function
function dydt=RHS(t,W,D,Nx,Ny,N,x,y,dx,dy,nu,k0,k1,a1,a2)
fprintf('Time: %g\n',t);
%reshape W into a matrix with 2*Ny rows and Nx-2 columns
wmat = reshape(W,[Nx-1 2*(Ny)])';
w = wmat(1:Ny,:); %displacements at interior nodes
v = wmat(Ny+1:2*Ny,:); %velocities at interior nodes

%initialize dydt
dydt = zeros(size([w;v]));
boundaryWxx = [ones(size(w,1),1) ones(size(w))];
boundaryWxxx = [ones(size(w,1),1) ones(size(w))];

dydt(1:Ny,1:Nx-1) = v; %place v values in the w_t locations

%pad v vector with zeros on all sides
v = [zeros(1,size(v,2)+2); [zeros(size(v,1),1), v, zeros(size(v,1),1)]; zeros(1,size(v,2)+2)];
v = [zeros(1,size(v,2)+2); [zeros(size(v,1),1), v, zeros(size(v,1),1)]; zeros(1,size(v,2)+2)];


gradarr=[zeros(size(w,1),1), w];
[gradx,grady] = grads(gradarr,Nx,Ny,dx,dy);
normgrad2 = trapz(y,trapz(x,gradx.*gradx + grady.*grady,2));

%pad w vector with zeros on all sides to calculate the grad w squared
w = [zeros(1,size(w,2)+2); [zeros(size(w,1),1), w, zeros(size(w,1),1)]; zeros(1,size(w,2)+2)];

w = [zeros(1,size(w,2)+2); [zeros(size(w,1),1), w, zeros(size(w,1),1)]; zeros(1,size(w,2)+2)];

%Clamped edge ghost points
w(:,1) = w(:,3);

%free edge (avoiding free-free corner area)
for j=2:Nx
    i = Ny+2;
    w(i+1,j) = -nu * (w(i,j-1) - 2 * w(i,j) + w(i,j+1)) / dx ^ 2 * dy ^ 2 - w(i-1,j) + 2 * w(i,j);
end

for j=2:Nx
    i = 3;
    w(i-1,j) = -nu * (w(i,j-1) - 2 * w(i,j) + w(i,j+1)) / dx ^ 2 * dy ^ 2 - w(i+1,j) + 2 * w(i,j);
end

for i=4:Ny+1
    j = Nx+1;
    w(i,j+1) = -nu * (w(i-1,j) - 2 * w(i,j) + w(i+1,j)) / dy ^ 2 * dx ^ 2 - w(i,j-1) + 2 * w(i,j);
end

%free edge, 2nd row (avoiding free-free corner area)

for j=3:Nx-1
    i = Ny+2;
    w(i + 2, j) = -(2 - nu) / dx ^ 2 * dy ^ 2 * (w(i + 1, j + 1) - w(i - 1, j + 1) - 2 * w(i + 1, j) + 2 * w(i - 1, j) + w(i + 1, j - 1) - w(i - 1, j - 1)) + 2 * w(i + 1, j) - 2 * w(i - 1, j) + w(i - 2, j);
end

for j=3:Nx-1
    i = 3;
    w(i - 2, j) = (2 - nu) / dx ^ 2 * dy ^ 2 * (w(i + 1, j + 1) - w(i - 1, j + 1) - 2 * w(i + 1, j) + 2 * w(i - 1, j) + w(i + 1, j - 1) - w(i - 1, j - 1)) + w(i + 2, j) - 2 * w(i + 1, j) + 2 * w(i - 1, j);
end

for i=5:Ny
    j = Nx+1;
    w(i,j+2) = -(2 - nu) * (w(i + 1, j + 1) - w(i + 1, j - 1) - 2 * w(i, j + 1) + 2 * w(i, j - 1) + w(i - 1, j + 1) - w(i - 1, j - 1)) / dy ^ 2 * dx ^ 2 + 2 * w(i, j + 1) - 2 * w(i, j - 1) + w(i, j - 2);
end

%Now we will apply the conditions for the corners.
%Labelling convention in this section is relative to the corner

%Top right corner
i = 3;
j = Nx + 1;
%E
w(i,j+1) = 2 * w(i,j) - w(i,j-1);

%EE
w(i, j + 2) = (-4 * (nu - 2) * (w(i, j) - w(i, j - 1) - w(i + 1, j + 1) / 2 + w(i + 1, j - 1) / 2) * dx ^ 2 + 4 * (w(i, j) - w(i, j - 1) + w(i, j - 2) / 4) * dy ^ 2) / dy ^ 2;

%N
w(i-1,j) = 2 * w(i,j) - w(i+1,j);

%NE
w(i-1,j+1) = w(i+1,j+1) - w(i+1,j-1) + w(i-1,j-1);

%NN
w(i-2, j) = (-4 * (nu - 2) * (w(i, j) - w(i + 1, j) - w(i - 1, j - 1) / 2 + w(i + 1, j - 1) / 2) * dy ^ 2 + 4 * dx ^ 2 * (w(i, j) - w(i + 1, j) + w(i + 2, j) / 4)) / dx ^ 2;

%NNW
w(i-2, j-1) = (2 * (w(i, j) - w(i + 1, j - 2) / 2 - w(i + 1, j) - w(i - 1, j - 1) + w(i + 1, j - 1) + w(i - 1, j - 2) / 2) * (nu - 2) * dy ^ 2 + 2 * dx ^ 2 * (w(i - 1, j - 1) - w(i + 1, j - 1) + w(i + 2, j - 1) / 2)) / dx ^ 2;

%SEE
w(i+1, j+2) = (2 * (w(i, j) - w(i, j - 1) - w(i + 1, j + 1) + w(i + 1, j - 1) + w(i + 2, j + 1) / 2 - w(i + 2, j - 1) / 2) * (nu - 2) * dx ^ 2 + dy ^ 2 * (w(i + 1, j - 2) + 2 * w(i + 1, j + 1) - 2 * w(i + 1, j - 1))) / dy ^ 2;


%Bottom right corner
i = Ny + 2;
j = Nx + 1;
%E
w(i,j+1) = 2 * w(i,j) - w(i,j-1);

%EE
w(i, j + 2) = (-4 * (nu - 2) * (w(i, j) - w(i, j - 1) - w(i - 1, j + 1) / 2 + w(i - 1, j - 1) / 2) * dx ^ 2 + 4 * (w(i, j) - w(i, j - 1) + w(i, j - 2) / 4) * dy ^ 2) / dy ^ 2;

%S
w(i+1,j) = 2 * w(i,j) - w(i-1,j);

%SE
w(i+1,j+1) = w(i-1,j+1) - w(i-1,j-1) + w(i+1,j-1);

%SS
w(i+2, j) = (-4 * (nu - 2) * (w(i, j) - w(i - 1, j) - w(i + 1, j - 1) / 2 + w(i - 1, j - 1) / 2) * dy ^ 2 + 4 * dx ^ 2 * (w(i, j) - w(i - 1, j) + w(i- 2, j) / 4)) / dx ^ 2;

%SSW
w(i+2, j-1) = (2 * (w(i, j) - w(i - 1, j - 2) / 2 - w(i - 1, j) - w(i + 1, j - 1) + w(i - 1, j - 1) + w(i + 1, j - 2) / 2) * (nu - 2) * dy ^ 2 + 2 * dx ^ 2 * (w(i + 1, j - 1) - w(i - 1, j - 1) + w(i - 2, j - 1) / 2)) / dx ^ 2;

%NEE
w(i-1, j+2) = (2 * (w(i, j) - w(i, j - 1) - w(i - 1, j + 1) + w(i - 1, j - 1) + w(i - 2, j + 1) / 2 - w(i - 2, j - 1) / 2) * (nu - 2) * dx ^ 2 + dy ^ 2 * (w(i - 1, j - 2) + 2 * w(i - 1, j + 1) - 2 * w(i - 1, j - 1))) / dy ^ 2;


%compute the RHS of the velocity equation
for j = 3:Nx+1 %which column of nodes
    for i=3:Ny+2 %which row of nodes

        %In vacuo conditions
        wxxxx = (1/dx^4) * (w(i-2,j) - 4*w(i-1,j) + 6*w(i,j)...
                - 4 *w(i+1,j) + w(i+2,j));
        
        wyyyy = (1/dy^4) * (w(i,j-2) - 4*w(i,j-1) + 6*w(i,j)...
                - 4 *w(i,j+1) + w(i,j+2));
        
        wxxyy = (1/(dx^2*dy^2)) * ( w(i+1,j+1) - 2* w(i+1,j)...
                + w(i+1,j-1) - 2* w(i,j+1) + 4 * w(i,j) -2*w(i,j-1)...
                +w(i-1,j+1) -2 * w(i-1,j)+w(i-1,j-1));
            
        %damping term
        wt = v(i,j);
        
        wtxx = (1/dx^2) * (v(i+1,j) - 2* v(i,j) + v(i-1,j));
        
        wtyy = (1/dy^2) * (v(i,j+1) - 2* v(i,j) + v(i,j-1));
        
        damping = - k0*(wt) + k1*(wtxx + wtyy);
        
        %spatial
        wx = (1/(2*dx))*(w(i,j+1)-w(i,j-1)); %w_x calculation
        wy = (1/(2*dy))*(w(i+1,j)-w(i-1,j)); %w_y calculation
        
        spatial = - a1 * wx - a2 * wy;
        
         dydt(Ny+(i-2),j-2) =  - 1 * D *(wxxxx + 2 * wxxyy + wyyyy)...
                                 + f(x(j-2),y(i-2),t)...
                                 + spatial + damping;  
    end
end
%dydt
dydt = reshape(dydt',[2*N 1]);

end
\end{lstlisting}

\newpage 

\nocite{*}

\printbibliography

\end{document}